\documentclass[letter, 11pt]{article}
\usepackage{amsfonts,latexsym,graphicx,amssymb,amsmath,url,amsthm}
\newtheorem{theorem}{Theorem}%[section]

\newtheorem{conjecture}{Conjecture}
\theoremstyle{definition}
\theoremstyle{remark}

\newcommand{\N}{\mathbb N}
\newcommand{\Z}{\mathbb Z}

\title{Golomb's conjecture on prime gaps}
\author{
Christian Elsholtz
}

\date{\vspace{-5ex}}
\begin{document}
\maketitle

\begin{abstract}
Question 10208b (1992) of %this Monthly 
the American Mathematical Monthly
asked: does there exist
an increasing  sequence 
$\{a_k\}$ of positive integers
and a constant $B > 0$
having the property that
$\{ a_k + n\}$ contains no more than $B$ primes for every integer $n$?
A positive answer to this question became known as Golomb's conjecture.
In this note we give a negative answer, making use of recent progress
in prime number theory.
\end{abstract}

\section{Introduction}
Solomon W.~Golomb 
\cite{Golomb-monthly, Golomb-proceedings-2010, Golomb-proceedings-2014}
repeatedly asked the following question, which was
Question 10208 (1992, page 266) of %this Monthly. 
the American Mathematical Monthly.

Let $1 < a_1 < a_2 < a_3 < \ldots$ be an increasing sequence of positive integers.
\begin{itemize}
\item[(a)] 
Is there such a sequence $\{a_k\}$ having the property that, for all integers 
$n$ (positive, negative, or zero), $\{a_k + n\}$ contains only finitely many
  primes?
\item[(b)] 
Is there such a sequence $\{a_k\}$ and a constant $B > 0$ 
having the property that 
$\{ a_k + n\}$ contains no more than $B$ primes for every integer $n$?
\end{itemize}

For Part a) a solution was given by Kevin Ford %(This Monthly, 1995, vol. 102, no. 4, p. 361--362). 
(American Mathematical Monthly, 1995, vol. 102, no. 4, p. 361--362). 
The answer is yes. Possible sequences are defined by $a_k = ((2k)!)! + k!$
or  $a_k=((2k)!)^3$.

For part b) Kevin Ford, Robert High, Gerry Myerson, and Solomom Golomb
 also pointed out that the following conjecture
would imply a negative answer.

\begin{conjecture}[Prime $k$-tuple conjecture]
Let $A=\{a_1, \ldots, a_k\}$ be a set of integers. If there is no prime $p$
such that the values $\{a_i \bmod p\}\
 (i=1, \ldots, k$) define a complete set of all $p$ residue classes,
 then $A$ is an admissible set.
For an admissible set there are infinitely many values $n$ such that
the $k$ values $a_i+n$ are simultaneously prime.
\end{conjecture}
The prime $k$-tuple conjecture is generally believed, but also 
considered to be hopelessly difficult. 
Ribenboim \cite{Ribenboim} called a positive answer to Part b) 
``Golomb's conjecture'', and Golomb \cite{Golomb-proceedings-2010}
himself took up this phrase.
Golomb discussed reasons to believe or disbelieve this conjecture. Let us
cite\footnote{Taken almost verbatim, with minor adaptions}
 his arguments in favour of this conjecture.

``The plausibility of Golomb's Conjecture arises from considering a sequence
$A = \{a_n\}$, all $n \geq 1$, with an extraordinarily fast growth rate. 
For example, if
$a_n =\left(\left(10^{10^{10^{10^{10^n}}}}\right)!\right)^3$,'' \ldots
%(compare this sequence to the examples given for  part a). 
``the expected number of primes in $A$, namely
$\sum_{n=1}^{\infty}\frac{1}{\ln a_n}$, is a tiny positive real
number, and this will be true for each translate sequence 
$A_{\tau} = \{a_n + \tau\}$, for all
$\tau\in \Z$. If we then take a huge value of $K$, e.g. 
$K_0 = 10^{10^{10^{10^{100}}}}$, for Golomb's
Conjecture to be false there must be {\emph{infinitely many}} values of 
$\tau$ with'' more than $K_0$ prime values of $a_n+\tau$.
%K0. (If there were only finitely many such values, 
%since CAP (τ) is always finite,
%there would be a largest value CAP (τ)MAX = L, and taking anyK >L, for this
%K we would have CAP (τ) <K for all τ ∈ Z.)
``Since far faster growth rates than the $A = \{a_n\}$ 
suggested here, and far larger
values than the $K_0$ suggested here, can be chosen, 
it takes considerable faith to
remain convinced of the truth of the prime $k$-tuples conjecture.''

The above definition of $a_n$ is a very thin subset of 
the sequence $a_k=((2k)!)^3$ mentioned before.

In a more recent update, Golomb \cite{Golomb-proceedings-2014} argues similarly,
referring to the fact that there is an uncountable number of sequences
increasing as fast as those above,
summarizing this as follows: ``Hence, it
should be extremely likely that there is at least one such sequence $A$ 
for which each of its translates $A+k$ contains no more than $B$ primes, 
for some finite bound $B$, where $B$ is allowed to be extremely large.''

Certainly not everybody would agree with these heuristic arguments, but it seemed
impossible to rigorously refute them.
Indeed,
the quite thin sequence of Fermat numbers $2^{2^i}+1, (i=0,1,\ldots) $ is generally
believed to contain only finitely many prime numbers.
The heuristic reason is as follows.
By the prime number theorem there are asymptotically $\frac{x}{\log x}+o(\frac{x}{\log x})$
primes $p \leq x$, and one may convert this information to argue
that a ``random integer'' $q$ is prime with probability $\frac{1}{\log q}$.
Hence, the ``expected number'' of primes in a ``random sequence'' $\{a_i\}, i\geq 1$
should be about $\sum_{i=1}^{\infty} \frac{1}{\log a_i}$.
Of course this model has its weaknesses, as for any algebraically 
defined sequence the 
members are not truly random. For example, for the Fermat numbers, 
$a_i=2^{2^i}+1$ 
some minor correction would be required,
taking into account that by construction all integers are odd, 
and are for $i\geq 1$ not divisible by three and so on.
But these corrections are minor, compared to the fact that 
$\sum_{i=0}^{\infty} \frac{1}{\log(2^{2^i}+1)}$
is convergent. 

Another weakness of Golomb's heuristic argument might be
to apply this argument also 
to many shifted copies of the same sequence. 
Even if the first sequence behaves according 
to the model, the shifted copies are certainly not independent.

In this note we will disprove Golomb's conjecture.
In order to do so, we make use of a \emph{weak} but proven version of 
the prime $k$-tuple conjecture.
In fact, there has been considerable
progress on gaps between two consecutive primes due to Zhang.
This is not quite what we need, as we need results on $k$-tuples, 
albeit weaker than the unproven $k$-tuple conjecture.
Below we state two results which were very recently proved by Maynard \cite{Maynard:2015, Maynard-dense}.
Both imply that Golomb's conjecture is wrong, thus giving a negative answer to part b) of
Problem 10208.

\begin{theorem}[Maynard \cite{Maynard:2015}, Theorem 1.2.]
 Let $m \in \N$. Let $r \in \N$ be sufficiently large depending on $m$, 
and let 
$A =\{a_1, a_2, \ldots, a_r\}$ be a set of $r$ distinct integers. 
Then there is a positive constant $c(m)$ such that:
\[
\frac{\# \left\{ \{h_1, \ldots, h_m\} \subseteq A : 
\text{for infinitely many $n$ all of $n + h_1, \ldots, n + h_m$ are prime }\right\} }
{\# \left\{ \{h_1, \ldots, h_m\} \subseteq A\right\}  }\geq c(m).\]
%(In other words, there are many 
%$m$-tuples which are good for our application).
\end{theorem}

Suppose there is an infinite sequence ${\cal A}$ satisfying Golomb's conjecture
with bound $B$. Choosing a large $r$, then
a positive proportion of all subsets of length $m=B+1$ of $\{a_1, \ldots,
a_r\}$ generate -
for infinitely many shifts - at least $m=B+1$ primes, 
contradicting the definition of the sequence ${\cal A}$.

Another theorem of Maynard has the same consequence, but has a different 
view point. We state only a special case of Maynard's much more general 
result.
\begin{theorem}[Maynard \cite{Maynard-dense}, Theorem 3.1]
There exists a constant $C$  such that if $r > C$
 and $A=\{a_1,\ldots , a_r\}$ is admissible, then
\[ \#\left\{n \leq  x : \#\{n + a_1, \ldots, n + a_r\} \text{ contains at least }
\frac{\log r}{C}  \text{ primes } \right\}\geq C_1 \frac{x}{(\log x)^r},\]
with some positive constant $C_1$.
\end{theorem}

In order to apply this Theorem to arbitrary integer sequences we first observe
that for this special case of Maynard's result one can remove the assumption that
the set $A$ is admissible.
Indeed, any set $A$ of $r\geq 25$ values $\{a_1, \ldots, a_r\}$ contains an 
admissible subset $A'\subset A$ of size $|A'|>\frac{r}{2\log r}$.
To see this, observe that modulo $2$ there exists a subset $A_1\subset A$
of size $|A_1|\geq |A|/2$ being only in one of the two residue classes.
Similarly, modulo $3$ there exists a subset  $A_2\subset A_1$
of size $|A_2|\geq |A_1|(1-\frac{1}{3})$ being only in two 
of the three residue classes. Iterating this modulo all primes $p \leq r$,
there exists a subset $A'\subset A$ of size
 $|A'|\geq |A|\prod_{p \leq r} (1-\frac{1}{p})$.
To simplify this last product over the primes we use 
an explicit estimate, 
namely formula (3.27) of Rosser and Schoenfeld \cite{RosserandSchoenfeld:1962}.
For all $r\geq 1$ the following holds:
\[\prod_{p \leq r}\left(1-\frac{1}{p}\right) > 
\frac{e^{-\gamma}}{\log r}\left (1-\frac{1}{(\log r)^2}\right), \]
where $\gamma = 0.577\ldots$ is the Euler-Mascheroni constant. 
If $r\geq 25$, then
\[\prod_{p \leq r}\left(1-\frac{1}{p}\right) > 
\frac{e^{-\gamma}}{\log r}\left (1-\frac{1}{(\log 25)^2}\right)\geq
\frac{1}{2\log r}. \]

We apply Maynard's theorem to the subsequence $A'\subset A$ with
$|A'|\geq \frac{r}{2\log r}$ elements.

For any given $B\geq 3$, choosing $r$ of size
 $\exp (2B C)$ suffices to guarantee that there are many values $n\leq x$
such that the $r$ values of the set $\{a_1+n, \ldots, a_r+n\}$
contain at least 
\[\frac{\log\left( \frac{r}{2\log r}\right)}{C}\geq 
\frac{\log\left( \frac{\exp(2BC)}{4BC}\right)}{C}\geq
2B- \frac{\log(4BC)}{C}>B
\]
primes, again contradicting the Golomb conjecture.

As a consequence, for any constant $B$ there are infinitely many values $n$
such that the sequence $\{2^{2^i}+n\}$ of translated Fermat numbers
contains at least $B+1$ primes.
However, we \emph{cannot} conclude that there is any fixed $n$ such that
the sequence $\{2^{2^i}+n\}$ contains infinitely many primes!

The author would like to thank Kaisa Matom\"aki for discussions on the 
subject, and the FIM at ETH Z\"urich, where this note was written up, 
for a very pleasant stay.

Christian Elsholtz\\
Institut f\"ur Analysis und Zahlentheorie\\
Technische Universit\"at Graz\\
Kopernikusgasse 24\\
A-8010 Graz\\ 
Austria.\\
{\it elsholtz@math.tugraz.at}
\ \\
\ \\
\noindent \textrm{MSC2010: primary: 11A41}

\end{document}